\documentclass{emsprocart}

\usepackage{mathrsfs,amssymb,amsfonts,amsmath}

\usepackage{graphics}
\usepackage[all]{xy}
\xyoption{curve}
\xyoption{import}
\xyoption{arc}
\xyoption{ps}
\UsePSspecials{dvips}
\usepackage{amsmath}
\usepackage{graphicx}
\usepackage[all]{xy}

\contact[carrillo@math.jussieu.fr]{Paulo Carrillo Rouse\\ Projet
  d'alg{\`e}bres d'op{\'e}rateurs\\ Universit{\'e} de Paris 7\\ 175, rue de
  Chevaleret\\ 
Paris, France}

\newtheorem{theorem}{Theorem}[section]

\newtheorem{lemma}[theorem]{Lemma}

\newtheorem{proposition}[theorem]{Proposition}

\theoremstyle{definition}

\newtheorem{definition}[theorem]{Definition}

\newtheorem{remark}[theorem]{Remark}

\newtheorem{examples}[theorem]{Examples}

\newcommand{\ci}{C^{\infty}}
\newcommand{\Cat}{\mathscr{C}}
\newcommand{\Dnc}{\mathscr{D}}
\newcommand{\gr}{\mathscr{G}}
\newcommand{\go}{\mathscr{G} ^{(0)}}

\newcommand{\Nb}{\mathscr{N}}
\newcommand{\sw}{\mathscr{S}}
\newcommand{\Uo}{\mathscr{U}}
\newcommand{\Vo}{\mathscr{V}}
\newcommand{\Rr}{\mathbf{R}}
\newcommand{\Nat}{\mathbf{N}}
\newcommand{\src}{\mathscr{S}_{r,c}}
\newcommand{\cg}{C_{c}^{\infty}(\gr)}
\newcommand{\cgo}{C_{c}^{\infty}(\go)}

\title[A Schwartz type algebra
  for the Tangent groupoid]{A Schwartz type algebra
  for the Tangent groupoid}

\author{Paulo Carrillo Rouse}

\begin{document}

\begin{abstract}
We construct an algebra of smooth functions over
the tangent groupoid associated to any Lie groupoid. This algebra is a
field of algebras over the closed interval $[0,1]$ which fiber at zero is the
algebra of Schwartz functions over the Lie algebroid, whereas any
fiber out of zero is the convolution algebra of the initial groupoid. 
Our motivation comes from index theory for Lie groupoids. In fact, our
construction gives an intermediate algebra between
the enveloping $C^*$-algebra and the convolution algebra
of compactly supported functions of the tangent groupoid; and it will
allows us, in a further work, to define other 
analytic index morphisms as a sort of ''deformations''. 
\end{abstract}

\begin{classification}

Primary 58-06, 19-06; Secondary 58H15, 19K56.

\end{classification}

\begin{keywords}

Lie groupoids, Tangent groupoid, K-theory, Index theory.

\end{keywords}

\maketitle

\section{Introduction}

The concept of groupoid is central in non commutative geometry. Groupoids
generalize the concepts of spaces, groups and 
equivalence relations. It is clear nowadays that 
groupoids are natural substitutes of singular
spaces. Many people have contributed to realizing this idea. We can
find for instance a groupoid-like treatment in Dixmier's works on
transformation groups, \cite{Dix}, or in Brown-Green-Rieffel's work on orbit
classification of relations, \cite{BGR}. 
In foliation theory, several models for
the leaf space of a foliation were realized using groupoids, mainly
by people like Haefliger (\cite{Ha}) and Wilkelnkemper (\cite{W}), for mention
some of them. There is also the case of Orbifolds, these can
be seen indeed as {\'e}tale groupoids, (see for example Moerdijk's paper 
\cite{Moer}). There are also some particular groupoid models for manifolds with
corners and conic manifoldss worked by people like Monthubert \cite{Mont},
Debord-Lescure-Nistor (\cite{DLN}) 
and Aastrup-Melo-Monthubert-Schrohe (\cite{AMMS}) for example.
 
The way we treat
''singular spaces'' in non commutative geometry is by associating to
them algebras.
In the case when the ''singular
space'' is represented by a Lie groupoid, we can, for instance,
consider the convolution algebra of differentiable functions with
compact support over the groupoid (see Connes or Paterson's books
\cite{Co2} and \cite{Pat}). This
last algebra plays the role of the algebra of smooth functions over
the ''singular space'' represented by the groupoid. From the
convolution algebra it is also possible to construct a $C^*$-algebra,
$C^*(\gr)$, that plays, in some sense, the role of the algebra of continuous
functions over the ''singular space''. The idea of associating
algebras in this sense can be traced back in works of Dixmier
(\cite{Dix}) for transformation groups, Connes (\cite{Co0}) for
foliations and Renault (\cite{Re}) for locally compact groupoids, for
mention some of them.
 
Using methods of Noncommutative Geometry, 
we would like to get invariants of this algebras, and hence, of the
spaces they represent. For that, 
Connes showed that many
groupoids and algebras associated to them appeared as `non commutative
analogues` of spaces to which many tools of geometry (and topology) such as
K-theory and Characteristic classes could be applied 
(\cite{Co1}, \cite{Co2}).
One classical way to obtain invariants in classical geometry
(topology), is through the index
theory in the sense of Atiyah-Singer. In the Lie groupoid case, there
is a Pseudodifferential calculus, developed by Connes (\cite{Co0}),
Monthubert-Pierrot (\cite{MP}) 
and Nistor-Weinstein-Xu (\cite{NWX}) in general. Some interesting
particular cases were treated in the groupoid-spirit by Melrose (\cite{Mel}),
Moroianu (\cite{Mor}) and others (see \cite{AMMS}).
Let $\gr \rightrightarrows \go$ be a Lie groupoid, there is an analytic
index morphism, (see \cite{MP}),
$$ 
ind_a:K^0(A^*\gr) \rightarrow K_0(C^*(\gr)),
$$
where $A\gr$ is the Lie algebroid of $\gr$. The
"$C^*$-index" $ind_a$ is a homotopy invariant of the  
$\gr$-pseudodifferential elliptic operators and has proved to be very
useful in very different situations, (see \cite{ANS}, \cite{CS}, \cite{DLN} for
example). One way to define the above index map is using the
Connes' tangent groupoid associated to $\gr$ as explained by Hilsum and
Skandalis in \cite{HS} or by Monthubert and Pierrot in \cite{MP}. The
tangent groupoid is a Lie groupoid
$$\gr^T \rightrightarrows \go \times [0,1]$$
with $\gr^T:= A\gr \times \{0\} \bigsqcup \gr \times (0,1]$ and the
groupoid structure is given by the groupoid structure of $A\gr$ at
$t=0$ and by the groupoid structure of $\gr$ for $t\neq 0$. One of the main
features about the tangent groupoid is that its $C^*-$algebra
$C^*(\gr^T)$ is a continuous field of $C^*-$algebras over the closed
interval $[0,1]$, with associated fiber algebras
\begin{center}
$C_0(A^*\gr)$ at $t=0$, and
\end{center}
\begin{center}
$C^*(\gr)$ for $t\neq 0$. 
\end{center}
In fact, it gives a $C^*-$algebraic
quantization of the Poisson manifold $A^*\gr$ (in the sense of
\cite{La}), and this is the main point why it allows to define the
index morphism as a sort of ''deformation''. Thus, the tangent
groupoid construction has been very useful in index theory 
(\cite{ANS}, \cite{DLN}, \cite{HS}) but also for other
purposes (\cite{La}, \cite{NWX}).

Now, to understand the purpose of the present work, let us first say
that the indices (in the sense of Atiyah-Singer-Connes) have not
necessarily  to be considered as elements in $K_0(C^*(\gr))$. Indeed,
it is possible to consider indices in $K_0(\ci_c (\gr))$. The
$\ci_c$-indices are more refined but they have several inconvenients
(see Alain Connes' book section $9.\beta$ 
for a discussion on this matter), nevertheless this
kind of indices have the great advantage that one can apply to them the
existent tools (such as pairings with cyclic cocycles or Chern-Connes
character) in order to obtain numerical invariants.

In this work we
begin a study of more refined indices. In particular we are looking
for indices between the $\ci_c$ and the $C^*$-levels; trying to keep
the advantages of both approaches 
(see \cite{Ca} for a more complete discussion). 
In the case of Lie groupoids this refinement could mean forget for a moment 
the powerful tools of the
theory of $C^*$-algebras and instead, working in a purely algebraic
and geometric level. In the present article, we construct 
an algebra of $\ci$ functions over
$\gr^T$, denoted by $\src (\gr^T)$. This algebra is also a 
field of algebras over the closed interval $[0,1]$, 
with associated fiber algebras,
\begin{center} 
$\sw(A\gr)$, at $t=0$, and
\end{center}
\begin{center}
$\ci_c(\gr)$ for $t\neq 0$,
\end{center}
where $\sw(A\gr)$ is the Schwartz algebra of the Lie
algebroid. Furthermore, we will have
\begin{equation}\label{incl}
\ci_c(\gr^T)\subset \src (\gr^T) \subset C^* (\gr^T),
\end{equation} 
as inclusions of algebras.
Let us explain in some words why we define an algebra over
the tangent groupoid such that in zero it is Schwartz: 
The "Schwartz algebras" have in general the good $K-$theory
groups. For example, we are interested in the symbols of $\gr$-PDO and
more precisely in their homotopy classes in $K$-theory, that is, we are
interested in the group $K^0(A^*\gr)=K_0(C_0(A^*\gr))$. Here it
would not be enough to take the $K-$theory of $\ci_c(A\gr)$ (see the 
example in \cite{Co2} p.142), however it is enough to consider the Schwartz
algebra $\sw (A^*\gr)$. Indeed, the Fourier transform shows 
that this last algebra is stable under holomorphic calculus
on $C_0(A^*\gr)$ and so it has the "good" $K$-theory, meaning that 
$K^0(A^*\gr)=K_0(\sw(A^*\gr))$. None of the inclusions in (\ref{incl})
is stable under holomorphic calculus, but that is precisely what we
wanted because our algebra $\src(\gr^T)$ have the remarkable property
that its evaluation at zero is stable under holomorphic calculus while
its evaluation at one (for example) is not.

The algebra $\src (\gr^T)$ is, as vector
space, a particular case of a more general construction that we do for
''Deformation to the normal cone manifolds'' from which the tangent groupoid
is a special case (see \cite{CC} and \cite{HS}). 
A deformation to the normal cone manifold (DNC for simplify) is a
manifold associated to an injective immersion $X\hookrightarrow M$
that is considered as a sort of blow up in differential geometry. The
construction of a DNC manifold has very nice functorial
properties (section \ref{dnc}) which we exploit to achieve our
construction. We think that our
construction 
could be used also for other purposes, 
for example, it seems that it could help to give more
understanding in quantization theory (see again \cite{CC}).

The article is organized as follows. 
In the second section we recall the basic facts about Lie groupoids. We
explain very briefly how to define the convolution algebra $\ci_c(\gr)$. 
In the third section we explain the "deformation to the normal cone"
construction associated to an injective immersion. Even if this could
be considered as classical material, we do it in some detail since we
will use in the sequel very explicit descriptions that we could not find
elsewhere. We also review some functorial properties associated to these
deformations.
A particular case
of this construction is the tangent groupoid associated to a Lie
groupoid. In the fourth section we start by constructing a vector
space $\src(\Dnc_{X}^{M})$ for any Deformation to the normal cone
manifold $\Dnc_{X}^{M}$; this space already exhibits the characteristic of
being a field of vector spaces over the closed interval $[0,1]$, such
that in zero we have a Schwartz space while out of zero we have
$\ci_c(M)$. We then define the algebra $\src
(\gr^T)$, the main result is precisely that the product is well
defined. 
The last section is devoted to motivate the construction of our
algebra by explaining in a few words some 
further developements that will immediately follow from this
work. All the results of the present work are part of the author's PHD thesis. 

\thanks{I want to thank my PHD advisor, Georges Skandalis, for
  all the ideas that he shared with me. I would also like 
to thank him for all the comments and remarks he made to the present
work. I would also like to thank the referee for the useful comments he made
for improving this paper.}

\section{Lie groupoids}

Let us recall what a groupoid is:

\begin{definition}
A $\it{groupoid}$ consists of the following data:
two sets $\gr$ and $\go$, and maps
\begin{itemize}
\item[$\cdot$]$s,r:\gr \rightarrow \go$ 
called the source and target map respectively,
\item[$\cdot$]$m:\gr^{(2)}\rightarrow \gr$ called the product map 
(where $\gr^{(2)}=\{ (\gamma,\eta)\in \gr \times \gr : s(\gamma)=r(\eta)\}$),
\item[$\cdot$]$u:\go \rightarrow \gr$ the unit map and 
\item[$\cdot$]$i:\gr \rightarrow \gr$
  the inverse map
\end{itemize}
such that, if we note $m(\gamma,\eta)=\gamma \cdot \eta$, $u(x)=x$ and 
$i(\gamma)=\gamma^{-1}$, we have 
\begin{itemize}
\item[1.]$\gamma \cdot (\eta \cdot \delta)=(\gamma \cdot \eta )\cdot \delta$, 
$\forall \gamma,\eta,\delta \in \gr$ when this is possible.
\item[2.]$\gamma \cdot x = \gamma$ and $x\cdot \eta =\eta$, $\forall
  \gamma,\eta \in \gr$ with $s(\gamma)=x$ and $r(\eta)=x$.
\item[3.]$\gamma \cdot \gamma^{-1} =u(r(\gamma))$ and 
$\gamma^{-1} \cdot \gamma =u(s(\gamma))$, $\forall \gamma \in \gr$.
\item[4.]$r(\gamma \cdot \eta) =r(\gamma)$ and $s(\gamma \cdot \eta) =s(\eta)$.
\end{itemize}
Generally, we denote a groupoid by $\gr \rightrightarrows \go $ where 
the parallel arrows are the source and target maps
and the other maps are given.
\end{definition}

Now, a Lie groupoid is a groupoid in which every set and map appearing
in the last definition 
is $\ci$ (possibly with borders), and the source and target maps are
submersions. For $A,B$ subsets of $\go$ we use the notation
$\gr_{A}^{B}$ for the subset $\{ \gamma \in \gr : s(\gamma) \in A,\, 
r(\gamma)\in B\}$.

All along this paper, $\gr \rightrightarrows \go $ is going to be a
Lie groupoid. We recall how to define an algebra structure in $\cg$ using
smooth Haar systems.
 
\begin{definition}
A $\it{smooth\, Haar\, system}$ over a Lie groupoid consists of a family of
measures $\mu_x$ in $\gr_x$ for each $x\in \go$ such that,
\begin{itemize}
\item for $\eta \in \gr_{x}^{y}$ we have the following compatibility
  condition:
$$\int_{\gr_x}f(\gamma)d\mu_x(\gamma)
=\int_{\gr_y}f(\gamma \circ \eta)d\mu_y(\gamma)$$
\item for each $f\in \cg$ the map
$$x\mapsto \int_{\gr_x}f(\gamma)d\mu_x(\gamma) $$ belongs to $\cgo$
\end{itemize}
\end{definition}

A Lie groupoid always posses a smooth Haar system. In fact, if we
fix a smooth (positive) section of the 1-density bundle associated to
the Lie algebroid we obtain a smooth Haar system 
in a canonical way. The advantage of using 1-densities is
that the measures are locally equivalent to the Lebesgue measure. 
We suppose for the rest of the
paper a given smooth Haar system given by 1-densities (for complete
details see \cite{Pat}). 
We can now define a convolution
product on $\cg$: Let $f,g\in \cg$, we set
$$(f*g)(\gamma)
=\int_{\gr_{s(\gamma)}}
f(\gamma \cdot \eta^{-1})g(\eta)d\mu_{s(\gamma)}(\eta)$$
This gives a well defined associative product. 
\begin{remark}
There
is a way to avoid the Haar system when one works with Lie groupoids,
using half densities (see Connes' book \cite{Co2}).
\end{remark}

\section{Deformation to the normal cone}\label{dnc}

Let $M$ be a $\ci$ manifold and $X\subset M$ be a $\ci$ submanifold. We denote
by $\Nb_{X}^{M}$ the normal bundle to $X$ in $M$, $\it{i.e.}$, 
$\Nb_{X}^{M}:= T_XM/TX$. We define the following set
\begin{align}
\Dnc_{X}^{M}:= \Nb_{X}^{M} \times {0} \bigsqcup M \times (0,1]. 
\end{align} 
The purpose of this section is to recall how to define a $\ci$-structure with
boundary in $\Dnc_{X}^{M}$. This is more or less classical, for example
it was extensively used in \cite{HS}. Here we are only going to
do a sketch.

Let us first consider the case where $M=\Rr^n$ 
and $X=\Rr^p \times \{ 0\}$ (where we
identify canonically $X=\Rr^p$). We denote by
$q=n-p$ and by $\Dnc_{p}^{n}$ for $\Dnc_{\Rr^p}^{\Rr^n}$ as above. In this case
we clearly have that $\Dnc_{p}^{n}=\Rr^p \times \Rr^q \times [0,1]$ (as a
set). Consider the 
bijection  
\begin{align}\label{psi}
\Psi: \Rr^p \times \Rr^q \times [0,1] \rightarrow
\Dnc_{p}^{n}
\end{align} 
given by 
$$
\Psi(x,\xi ,t) = \left\{ 
\begin{array}{cc}
(x,\xi ,0) &\mbox{ if } t=0 \\

(x,t\xi ,t) &\mbox{ if } t>0
\end{array}\right.
$$
which inverse is given explicitly by 
$$
\Psi^{-1}(x,\xi ,t) = \left\{ 
\begin{array}{cc}
(x,\xi ,0) &\mbox{ if } t=0 \\
(x,\frac{1}{t}\xi ,t) &\mbox{ if } t>0
\end{array}\right.
$$
We can consider the $\ci$-structure with border on $\Dnc_{p}^{n}$
induced by this bijection.  

In the general case. Let 
$(\Uo,\phi)$ be a local chart in $M$ and suppose it is an $X$-slice,
so that it
satisfies
\begin{itemize}
\item[1)]$\phi : \Uo \stackrel{\cong}{\rightarrow} U \subset \Rr^p\times \Rr^q$
\item[2)]If $\Uo \cap X =\Vo$, $\Vo=\phi^{-1}( U \cap \Rr^p \times \{ 0\}
  )$ (we note $V=U \cap \Rr^p \times \{ 0\}$)
\end{itemize}
With this notation we have that $\Dnc_{V}^{U}\subset \Dnc_{p}^{n}$ is an
open subset. We may define a function $$
\tilde{\phi}:\Dnc_{\Vo}^{\Uo} \rightarrow \Dnc_{V}^{U} 
$$ in the following way: For $x\in \Vo$ we have $\phi (x)\in \Rr^p
\times \{0\}$. If we write 
$\phi(x)=(\phi_1(x),0)$, then 
$$ \phi_1 :\Vo \rightarrow V \subset \Rr^p$$ 
is a diffeomorphism, where $V=U\cap (\Rr^p \times \{0\})$. We set 
$\tilde{\phi}(v,\xi ,0)= (\phi_1 (v),d_N\phi_v (\xi ),0)$ and 
$\tilde{\phi}(u,t)= (\phi (u),t)$ 
for $t\neq 0$. Here 
$d_N\phi_v: \Nb_v \rightarrow \Rr^q$ is the normal component of the
 derivate $d\phi_v$ for $v\in \Vo$. It is clear that $\tilde{\phi}$ is
 also a  bijection (in particular it induces a $C^{\infty}$ structure 
with border over $\Dnc_{\Vo}^{\Uo}$).
 
Let us define, with the same notations as above, the following set 
$$\Omega_{V}^{U}=\{(x,\xi,t)\in 
\Rr^p \times \Rr^q \times [0,1]: (x,t\cdot \xi)\in U \}. $$ 
which is an open subset of $\Rr^p \times \Rr^q \times [0,1]$ and thus
a $\ci$ manifold (with border). 
It is immediate that $\Dnc_{V}^{U}$ is diffeomorphic to $\Omega_{V}^{U}$
through the restriction of $ \Psi$, used in (\ref{psi}). 
Now we consider an atlas 
$ \{ (\Uo_{\alpha},\phi_{\alpha}) \}_{\alpha \in \Delta}$ of $M$
 consisting of $X-$slices. It is clear that 
\begin{equation}\label{atlasdcn} 
\Dnc_{X}^{M}= \cup_{\alpha \in
 \Delta}\Dnc_{\Vo_{\alpha}}^{\Uo_{\alpha}}
\end{equation}
and if we take $\Dnc_{\Vo_{\alpha}}^{\Uo_{\alpha}} 
\stackrel{\varphi_{\alpha}}{\rightarrow } 
\Omega_{V_{\alpha}}^{U_{\alpha}}$ defined as the composition
$$\Dnc_{\Vo_{\alpha}}^{\Uo_{\alpha}} \stackrel{\phi_{\alpha}}{\rightarrow}
\Dnc_{V_{\alpha}}^{U_{\alpha}}
\stackrel{\Psi_{\alpha}^{-1}}{\rightarrow} 
\Omega_{V_{\alpha}}^{U_{\alpha}} $$
then we obtain the following result.

\begin{proposition}\label{atlas}
$ \{ (\Dnc_{\Vo_{\alpha}}^{\Uo_{\alpha}},\varphi_{\alpha})
  \} _{\alpha \in \Delta }$ is a $\ci$ atlas with border over
  $\Dnc_{X}^{M}$.
\end{proposition}

In fact the proposition can be proved directly from the following
elementary lemma

\begin{lemma}\label{ellema2}
Let $F:U\rightarrow U'$ a $\ci$ diffeomorphism where 
$U\subset \Rr^p \times \Rr^q$ and $U'\subset \Rr^{p} \times \Rr^{q}$ are
open subsets. We write $F=(F_1,F_2)$ and we suppose that
$F_2(x,0)=0$. Then the function 
$\tilde{F}:\Omega_{V}^{U} \rightarrow \Omega_{V'}^{U'}$ defined by
$$
\tilde{F}(x,\xi ,t) = \left\{ 
\begin{array}{cc}
(F_1(x,0),\frac{\partial F_2}{\partial \xi}(x,0) \cdot \xi,0) 
&\mbox{ if } t=0 \\
(F_1(x,t\xi),\frac{1}{t}F_2(x,t\xi),t) &\mbox{ if } t>0
\end{array}\right.
$$
is a $\ci$ map.
\end{lemma}

\begin{proof}
Since the result will hold if and only if it is true in each
coordinate, it is enough to prove that if we have $F:U \rightarrow  \Rr$ a
$C^{\infty}$ map with $F(x,0)=0$, then the map $\tilde{F}:\Omega_{V}^{U} 
\rightarrow \Rr $ given by
$$
\tilde{F}(x,\xi ,t) = \left\{ 
\begin{array}{cc}
\frac{\partial F}{\partial \xi}(x,0) \cdot \xi &\mbox{ if } t=0 \\
\frac{1}{t}F(x,t\xi) &\mbox{ if } t>0
\end{array}\right.
$$
is a $C^{\infty}$ map. For that, we write
\[ F(x,\xi)= \frac{\partial F}{\partial \xi}(x,0) \cdot \xi +
h(x,\xi)\cdot \xi\]
with $h:U \rightarrow \Rr^q$ a
$C^{\infty}$ map such that $h(x,0)=0$. Then 
\[ \frac{1}{t}F(x,t\xi)= \frac{\partial F}{\partial \xi}(x,0) \cdot \xi +
h(x,t\xi)\cdot \xi \]
from which we immediately get the result.
\end{proof}

\begin{definition}[DNC]
Let $X\subset M$ be as above. The set
$\Dnc_{X}^{M}$ provided with the  $C^{\infty}$ structure with border
induced by the atlas described in the last proposition is called
$\it{`The\, deformation\, to\, normal\, cone\, associated\, to\,}$   
$X\subset M$`. We will often write DNC instead of
Deformation to the normal cone. 
\end{definition}

\begin{remark}
Following the same steps, 
it is possible to define a deformation to the normal
cone associated to an injective immersion $X\hookrightarrow M$.
\end{remark}

\begin{examples}\label{exdcn}
Let us mention some basic examples of DCN manifolds 
$\Dnc_{X}^{M}$:
\begin{itemize}
\item[1.]Consider the case when $X=\emptyset$. We have that 
$\Dnc_{\emptyset}^{M}= M\times (0,1]$ with the usual $\ci$ structure
on $M\times (0,1]$. We used this fact implicitly for cover
$\Dnc_{X}^{M}$ as in (\ref{atlasdcn}). 
\item[2.]Consider the case when $X\subset M$ is an open subset. Then
  we do not have any deformation at zero and we immediately see by
  definition that $\Dnc_{X}^{M}$ is just the open subset of $M\times
  [0,1]$ consisting in the union of $X\times[0,1]$ and $M\times (0,1]$.
\end{itemize}
\end{examples}

The most important feature about the DNC construction is that it is in
some sense functorial. More explicitly, let $(M,X)$ 
and $(M',X')$ be $\ci$-couples as above and let
 $F:(M,X)\rightarrow (M',X')$
be a couple morphism, i.e., a $\ci$ map   
$F:M\rightarrow M'$, with $F(X)\subset X'$. We define 
$ \Dnc(F): \Dnc_{X}^{M} \rightarrow \Dnc_{X'}^{M'} $ by the following formulas:\\

$\Dnc(F) (x,\xi ,0)= (F(x),d_NF_x (\xi),0)$ and\\

$\Dnc(F) (m ,t)= (F(m),t)$ for $t\neq 0$,\\

where $d_NF_x$ is by definition the map
\[  (\Nb_{X}^{M})_x 
\stackrel{d_NF_x}{\longrightarrow}  (\Nb_{X'}^{M'})_{F(x)} \]
induced by $ T_xM 
\stackrel{dF_x}{\longrightarrow}  T_{F(x)}M'$.

We have the following proposition, which is also an immediate
consequence of the lemma above.
\begin{proposition}
The map $\Dnc(F):\Dnc_{X}^{M} \rightarrow \Dnc_{X'}^{M'}$ is $\ci$.
\end{proposition}

\begin{remark}
If we consider the category $\Cat_{2}^{\infty}$ of  $\ci$ pairs given by
a $\ci$ manifold and a $\ci$ submanifold and pair morphisms as above,
we can reformulate the proposition and say that we have a functor
$$\Dnc : \Cat_{2}^{\infty} \rightarrow \Cat^{\infty}$$ where
$\Cat^{\infty}$ denote the category of  $\ci$ manifolds with border.
\end{remark}

\subsection{The tangent groupoid}

\begin{definition}[Tangent groupoid]
Let $\gr \rightrightarrows \go $ be a Lie groupoid. \linebreak 
$\it{The\, tangent\,
groupoid}$ associated to $\gr$ is the groupoid that has $\Dnc_{\go}^{\gr}
$ as the set of arrows and  $\go \times [0,1]$ as the units, with:
\begin{itemize}
\item[$\cdot$] $s^T(x,\eta ,0) =(x,0)$ and $r^T(x,\eta ,0) =(x,0)$ at $t=0$.
\item[$\cdot$] $s^T(\gamma,t) =(s(\gamma),t)$ and $r^T(\gamma,t)
  =(r(\gamma),t)$ at $t\neq0$.
\item[$\cdot$] The product is given by
  $m^T((x,\eta,0),(x,\xi,0))=(x,\eta +\xi ,0)$ and \linebreak $m^T((\gamma,t), 
  (\beta ,t))= (m(\gamma,\beta) , t)$ if $t\neq 0 $ and 
if $r(\beta)=s(\gamma)$.
\item[$\cdot$] The unit map $u^T:\go \rightarrow \gr^T$ is given by
 $u^T(x,0)=(x,0)$ and $u^T(x,t)=(u(x),t)$ for $t\neq 0$.
\end{itemize}
We denote $\gr^{T}:= \Dnc_{\go}^{\gr}$.
\end{definition} 
As we have seen above $\gr^{T}$ can be considered as a $\ci$ manifold with
border. As a consequence of the functoriality of the DNC construction
we can show that the tangent groupoid is in fact a Lie
groupoid. Indeed, it is easy to check that if we identify in a
canonical way $\Dnc_{\go}^{\gr^{(2)}}$ with $(\gr^T)^{(2)}$, then 
$$ m^T=\Dnc(m),\, s^T=\Dnc(s), \,  r^T=\Dnc(r),\,  u^T=\Dnc(u)$$
where we are considering the following pair morphisms:
\begin{align}  
m:((\gr)^{(2)},\go)\rightarrow (\gr,\go ), \nonumber
\\
s,r:(\gr ,\go) \rightarrow (\go,\go),\nonumber 
\\
u:(\go,\go)\rightarrow (\gr,\go ).\nonumber
\end{align}
Finally, if $\{ \mu_x\}$ is a smooth Haar system on $\gr$, then, setting
\begin{itemize}
\item $\mu_{(x,0)}:=\mu_x$ at $(\gr^T)_{(x,0)}=T_x\gr_x$ and
\item $\mu_{(x,t)}:=t^{-q}\cdot \mu_x $ at $(\gr^T)_{(x,t)}=\gr_x$ for
  $t\neq 0$, where $q=dim\, \gr_x$, 
\end{itemize}
one obtains a smooth Haar system for the Tangent groupoid 
(details may be found in \cite{Pat}).

\begin{examples}\label{exgt}
We finish this section with some interesting examples of groupoids and their
tangent groupoids.
\item[$(i)$]$\it{The\, tangent\, groupoid\, of\, a\, group}$. 
Let $G$ be a Lie group
  considered as a Lie groupoid, $\gr:= G\rightrightarrows \{ e\}$. In this
  case the normal bundle to the inclusion $\{ e\}\hookrightarrow G$ is
  of course identified with the Lie algebra of the Group. Hence, the
  tangent groupoid is a deformation of the group in its Lie algebra:
\begin{center}
$\gr^T=\mathfrak{g}\times \{ 0\}\bigsqcup G\times (0,1]$.
\end{center}
\item[$(ii)$]$\it{The\, tangent\, groupoid\, of\, a\,
    smooth\, vector\, bundle}$. Let $E\stackrel{p}{\rightarrow} X$ 
be a smooth vector bundle over a $\ci$ manifold $X$ (connexe). We can
consider the Lie groupoid $E\rightrightarrows X$ induced by the vector
structure of the fibers, $\it{i.e.}$, $s(\xi)=p(\xi)=r(\xi)$ and the
composition is given by the vector sum $\xi \circ \eta =\xi+\eta$. In
this case the normal vector bundle associated to the zero section can
be identified to $E$ itself. Hence, as a set the tangent groupoid is
$E\times [0,1]$ but the $\ci$-structure at zero 
is given locally as in (\ref{psi}).  
\item[$(iii)$]$\it{The\, tangent\, groupoid\, of\, a\,
    \ci-manifold}$. Let $M$ a $\ci$-manifold. We can consider the
  product groupoid $\gr_M:=M\times M\rightrightarrows M$. The tangent
  groupoid in this case takes the following form
\begin{center}
$\gr^T_M=TM \times \{ 0\} \bigsqcup M\times M\times (0,1]$.
\end{center} 
This is called the tangent groupoid to $M$ and it was introduced by
Connes for giving a very conceptual proof of the Atiyah-Singer index
theorem (see \cite{Co2} and \cite{DLN}).
\end{examples}

\section{An algebra for the Tangent groupoid}\label{srcdnc}
In this section we will show how to construct an algebra for the
tangent groupoid which consist of $\ci$ functions that satisfy a rapid
decay condition at zero while out of zero they satisfy a compact
support condition. This algebra is the main construction in
this work.
\subsection{Schwartz type spaces for Deformation to the normal cone manifolds}
Our algebra for the Tangent groupoid will be a particular
case of a construction associated to any deformation to the normal
cone. We start by defining a space for DNCs associated to
open subsets of $\Rr^p \times \Rr^q$. 

\begin{definition}\label{ladef}
Let $p, q\in \Nat$ and $U \subset \Rr^p \times \Rr^q$ 
an open subset, and let $V=U\cap (\Rr^p \times \{ 0\})$.
\begin{itemize}
\item[(1)]Let $K\subset U \times [0,1]$ be a compact
  subset. We say that $K$ is a conic compact subset of $U \times [0,1]$
relative to $V$ if
\[ K_0=K\cap (U \times \{ 0\}) \subset V\]
\item[(2)]Let $g \in \ci
  (\Omega_{V}^{U})$. We say that
   $f$ has compact conic support $K$, if there exists a conic
  compact $K$
 of $U \times [0,1]$ relative to $V$ such that if $t\neq 0$ and 
$(x, t\xi ,t) \notin K$ then $g(x, \xi ,t)=0$.
\item[(3)]We denote by $\src (\Omega_{V}^{U})$ 
the set of functions
$g\in \ci (\Omega_{V}^{U})$ 
that have compact conic support and that satisfy the following condition:
\begin{itemize}
\item[$(s_1$)]$\forall$ $k,m\in \Nat$, $l\in \Nat^p$
and $\alpha \in \Nat^q$ it exists $C_{(k,m,l,\alpha)} >0$ such that
\[ (1+\| \xi \|^2)^k \| \partial_{x}^{l}\partial_{\xi}^{\alpha}
\partial_{t}^{m}g(x,\xi ,t) \| \leq C_{(k,m,l,\alpha)}   \]
\end{itemize}
\end{itemize}
\end{definition}
Now, the spaces $\src (\Omega_{V}^{U})$ are invariant under
diffeomorphisms. More precisely if $F:U\rightarrow U'$ is a $\ci$
diffeomorphism as in lemma \ref{ellema2} then we can prove the next result.
\begin{proposition}\label{gtilde}
Let $g\in \src (\Omega_{V'}^{U'})$, then 
$\tilde{g}:= g\circ \tilde{F} \in \src (\Omega_{V}^{U})$.
\end{proposition}

\begin{proof}
The first observation is that $\tilde{g} \in \ci (\Omega_{V}^{U})$,
thanks to lemma \ref{ellema2}. Let us check that it has compact conic
support. For that, let $K'\subset U'\times [0,1]$ the conic compact
support of $g$. We let 
$$K=(F^{-1}\times id_{[0,1]})\subset U\times [0,1],$$
which is a conic compact subset of $U\times[0,1] $ relative to $V$,
and it is immediate by definition that $\tilde{g}(x,\xi,t)=0$ if
$t\neq 0$ and $(x,t\cdot \xi,t)\notin K$, that is, $\tilde{g} $ has
compact conic support $K$.
 
We now check the rapid decay
property $(s_1)$: For simplify the proof we first introduce 
some useful notation. Writing $F=(F_1,F_2) $ as in the lemma
\ref{ellema2}, we denote
$F_1(x,\xi)=(A_1(x,\xi),...,A_p(x,\xi))$ and 
$F_2(x,\xi)=(B_1(x,\xi),...,B_q(x,\xi))$. We denote also
$w=w(x,\xi,t)=(A_1(x,t\xi),...,A_p(x,t\xi))$ and
$\eta=\eta (x,\xi,t)=(\tilde{B_1}(x,\xi,t),...,\tilde{B_q}(x,\xi,t)
)$ where $\tilde{B_j} $ is also as above, $\it{i.e.}$, 
$$
\tilde{B_j}(x,\xi ,t) = \left\{ 
\begin{array}{cc}
\frac{\partial B_j}{\partial \xi}(x,0)\cdot \xi &\mbox{ if } t=0 \\
\\
\frac{1}{t}B_j(x,t\xi) &\mbox{ if } t\neq 0\\
\end{array}\right.
$$
In particular by definition we have
$\tilde{F}(x,\xi,t)=(w,\eta,t)$. We also write $z=(x,\xi,t)$ 
and $u=(\omega,\eta,t)$. Hence, what we would like is to find bounds
for expressions of the following type
$$\| \xi \|^k\|\partial_{z}^{\alpha}\tilde{g}(z) \|,$$
for arbitrary $k\in \Nat$ and $\alpha \in \Nat^p \times \Nat^q \times
\Nat$. A simple calculation shows that the derivates 
$\partial_{z}^{\alpha}\tilde{g}(z)$ are of the following form
$$\partial_{z}^{\alpha}\tilde{g}(z)=
\sum_{|\beta| \leq |\alpha|} P_{\beta}(z)\partial_{u}^{\beta}g(u)$$
where $P_{\beta}(z)$ is a finite sum of products of the form
$$\partial_{z}^{\gamma}\omega_i(z)\cdot \partial_{z}^{\delta}\eta_j(z).$$
We are only interested in see what happens in the set 
$K_{\Omega}:=\{z=(x,\xi,t)\in \Omega:(x,t\cdot \xi,t)\in K \}$ since
out of this set we have 
that $g$ and all its derivates vanish ($(x,t\xi,t) \in K$ iff 
$(w,t\eta,t) \in K'$). For a point $z=(x,\xi,t)\in
K_{\Omega}$ we have that $(x,t\cdot \xi)$ is in a compact set and then
it follows that the expressions 
$$\| \partial_{z}^{\gamma}\omega_i(z) \|$$ are bounded in $K_{\Omega}$.
For the expressions $\| \partial_{z}^{\delta}\eta_j(z)\|$, we
proceed first by developing as in lemma \ref{ellema2}, that is, 
$$\eta_j(x,\xi,t)
=\left( \frac{\partial B_j}{\partial \xi}(x,0) \cdot \xi +
h^j(x,t\xi)\right) \cdot \xi $$
Now, since we are only considering points in $K_{\Omega}$, it is
immediate that we can find constants $C_j>0$ such that
$$\| \partial_{z}^{\delta}\eta_j(z)\|\leq C_j \cdot \| \xi \|^{m_{\delta}}.$$
In the same way (remember $F$ is a diffeomorphism) we can have
constants $C_i>0$ such that 
$$\| \xi_i (\omega,\eta,t)\| \leq C_i \cdot \| \eta \|.$$
Putting all together, and using the property $(s_1)$ for $g$, we
get bounds $C>0$ such that
$$\| \xi \|^k\|\partial_{z}^{\alpha}\tilde{g}(z) \| \leq C,$$
and this concludes the proof. 
\end{proof}

\begin{remark}
We can resume the last invariance result as follows: If $(\Uo,\Vo)$ is
a $\ci$ pair diffeomorphic to $(U,V)$ with $U\subset E$, an open
subset of a vector space $E$, and $V=U\cap E$, then
$\src(\Dnc_{\Vo}^{\Uo})$ is well defined and does not depend on the
pair diffeomorphism.
\end{remark}
With the last compatibility result in hand we are ready to give the
main definition in this work.

\begin{definition}\label{src}
Let $g \in \ci (\Dnc_{X}^{M}) $.
\begin{itemize}
\item[(a)]We say that $g$ has compact conic support $K$, 
if there exists a compact subset
 $K\subset M \times [0,1]$ with $K_0:=K\cap (M\times \{ 0\}) \subset X$ (conic
 compact relative to $X$) such that if $t\neq 0$ and 
$(m,t) \notin K$ then $g(m,t)=0$.
\item[(b)]We say that $g$ is rapidly decaying at zero if for every
$(\Uo,\phi)$  $X$-slice chart
and for every $\chi \in \ci_c(\Uo \times [0,1])$, the map 
$g_{\chi}\in
\ci(\Omega_{V}^{U})$
given by
\[ g_{\chi}(x,\xi ,t)= (g\circ \varphi^{-1})(x,\xi ,t) 
\cdot (\chi \circ p \circ \varphi^{-1})(x,\xi ,t) \]
is in  $\src (\Omega_{V}^{U})$, where $p$ is the projection 
$p:\Dnc_{X}^{M} \rightarrow M\times
[0,1]$ given by $(x,\xi,0)\mapsto (x,0)$, and 
$(m,t)\mapsto (m,t)$ for $t\neq 0$.
\end{itemize}
Finally, we denote by $\src (\Dnc_{X}^{M})$ the set of functions 
$g\in \ci(\Dnc_{X}^{M})$ that are rapidly decaying at zero 
with compact conic support.
\end{definition}

\begin{remark}
\begin{itemize}
\item[(a)] By definition of $\src (\Dnc_{X}^{M})$ we see that 
$\ci_c (\Dnc_{X}^{M})$ is contained as a vector subspace. 
\item[(b)] It is clear that $\ci_c(M\times(0,1])$ 
can be considered as a subspace 
of $\src (\Dnc_{X}^{M})$ by extending by zero the functions at $\Nb_{X}^{M}$.
\end{itemize}
\end{remark}

Following the lines of the last remark we are going to precise a
possible decomposition of our space $\src (\Dnc_{X}^{M})$ that will be
very useful in the sequel. 
Let 
$\{ (\Uo_{\alpha},\phi_{\alpha}) \}_{\alpha \in \Delta}$ a family of
$X-$slices covering $X$. Consider the open cover of $M\times [0,1]$
consisting in $\{ (\Uo_{\alpha}\times [0,1],\phi_{\alpha}) \}_{\alpha \in \Delta}$
union with $M\times (0,1]$. We can take a partition of the unity
subordinated to the last cover, 
\[ 
\{ \chi_{\alpha},\lambda\}_{\alpha \in
  \Delta}. 
\]
That is, we have the following properties:
\begin{itemize}
\item[$\cdot$]$0\leq \chi_{\alpha},\lambda \leq 1$
\item[$\cdot$]$\text{supp}\, \chi_{\alpha} \subset \Uo_{\alpha}\times [0,1]$ and 
$\text{supp}\, \lambda \subset M\times(0,1]$.
\item[$\cdot$]
$\sum_{\alpha}\chi_{\alpha}+\sum \lambda =1$
\end{itemize}
Let $f\in \src (\Dnc_{X}^{M})$, we denote
$$f_{\alpha}:=f|_{\Dnc_{\Vo_{\alpha}}^{\Uo_{\alpha}}}
\cdot (\chi_{\alpha}\circ p) \in \ci (\Dnc_{\Vo_{\alpha}}^{\Uo_{\alpha}})$$ and
$$f_{\lambda}:=f|_{M\times (0,1]} 
\cdot (\lambda \circ p) \in \ci (M\times
(0,1]),$$
then we obtain the following decomposition:
\[ f=\sum_{\alpha}f_{\alpha}+
  f_{\lambda} \]
Now, since $f$ is conic compactly supported we can suppose, without
lost of generality, that
\begin{itemize}
\item[$\cdot$]$f_{\lambda} \in \ci_c(M\times (0 ,1])$, and
\item[$\cdot$]that
  $\chi_{\alpha}$ is compactly supported in $\Uo_{\alpha} \times [0,1]$.
\end{itemize}
What we conclude of all this, is that we can decompose our space 
$\src (\Dnc_{X}^{M})$ as follows
\begin{align}\label{decomposicion}
\src (\Dnc_{X}^{M}) = \sum_{\alpha \in \Lambda} \src (\Dnc_{\Vo_{\alpha}}^{\Uo_{\alpha}}) 
+ \ci_c (M\times (0,1]).
\end{align}

As we mentioned in the introduction, we want to see the space $\src
(\Dnc_{X}^{M})$ as a field of vector spaces over the interval $[0,1]$, 
where at zero we talked
about Schwartz spaces. In our case we are interested in Schwartz
functions on the vector bundle $\Nb_{X}^{M}$.
Let us first recall the notion of the Schwartz space associated to a
vector bundle.

\begin{definition} 
Let $(E,p,X)$ be a smooth vector bundle over a $\ci$ manifold $X$. We
define the Schwartz space $\sw (E)$ as the set of $\ci$ functions
$g \in \ci (E)$ such that $g$ is a Schwartz function at each fiber (uniformly)
and $g$ has compact support in the direction of $X$, $\it{i.e.}$, if
there exists a compact subset $K\subset X$ such that $g(E_x)=0$ for
$x\notin K$.
\end{definition}
The vector space $\sw(E)$ is an associative algebra with the product
given as follows: for $f,g\in \sw(E)$, we put
\begin{align}
(f*g)(\xi)=\int_{E_{p(\xi)}}f(\xi-\eta)g(\eta)d\mu_{p(\xi)}(\eta),
\end{align}
where $\mu_{\xi}$ is a smooth Haar system of the Lie groupoid
$E\rightrightarrows X$. A classical Fourier argument can be applied to
show that the last algebra is isomorphic to $(\sw(E^*), \cdot)$
(punctual product). In particular this implies that $K_0(\sw(E))\cong
K^0(E^*)$.

In the case we are interested, we have a couple $(M,X)$ and a vector
bundle associated to it, that is, the normal bundle over $X$,
$\Nb_{X}^{M}$. The reason why we gave the last definition is because
we get evaluation linear maps 
\begin{align}\label{e0}
e_0:\src (\Dnc_{X}^{M}) \rightarrow \sw (\Nb_{X}^{M}),
\end{align} 
and
\begin{align}\label{e1}
e_t:\src (\Dnc_{X}^{M}) \rightarrow \ci_c(M)
\end{align} 
for $t\neq0$.
Consequently, we have that the vector space $\src (\Dnc_{X}^{M})$ is a
field of vector spaces over the closed interval $[0,1]$, which fibers
spaces are: $\sw (\Nb_{X}^{M})$ at $t=0$ and $\ci_c(M)$ for $t\neq 0$.

\begin{examples}
Let us finish this subsection by giving the examples of spaces 
$\src (\Dnc_{X}^{M})$ corresponding to the DCN manifolds seen at
\ref{exdcn} above.
\begin{itemize}
\item[1.]For $X=\emptyset$, we have that $\src
  (\Dnc_{\emptyset}^{M})\cong \ci_c(M\times (0,1])$.
\item[2.]For $X\subset M$ an open subset we have that $\src
  (\Dnc_{X}^{M})\cong \ci_c(W)$ where $W \subset M\times[0,1]$ is the
  open subset consisting of the union of $X\times [0,1]$ and $M\times (0,1]$.
\end{itemize}
\end{examples}

\subsection{Schwartz type algebra for the Tangent groupoid}

In this section we define an algebra structure on $\src (\gr^T)$. 
We start by defining a function 
$m_{r,c}: \src (\Dnc_{\go}^{\gr^{(2)}})
\rightarrow \src (\Dnc_{\gr^{(0)}}^{\gr})$ by the following formulas:

For $F\in \src (\Dnc_{\go}^{\gr^{(2)}})$, we let
\[ m_{r,c}(F)(x,\xi ,0)= \int_{T_x\gr_x}F(x, \xi -\eta ,
\eta,0)
d\mu_x(\eta) \] 
and
\[ m_{r,c}(F)(\gamma,t)= \int_{\gr_{s(\gamma)}}F(\gamma \circ \delta^{-1},
\delta ,t) t^{-q}d\mu_{s(\gamma)}(\delta) \]
If we canonically identify $\Dnc_{\go}^{\gr^{(2)}}$ with
$(\gr^T)^{(2)}$, the map above is nothing else that the integration
along the fibers of $m^T:(\gr^T)^{(2)}\rightarrow \gr^T$.
We have the following proposition:

\begin{proposition}\label{elproducto}
$m_{r,c}: \src ((\gr^T)^{(2)})
\rightarrow \src (\gr^T)$ is a well defined linear map.
\end{proposition}

The interesting part of the proposition is that the map is well
defined since it will evidently be linear. Let us suppose for the
moment that the last proposition is true.
Under this assumption, we will define the product in $\src (\gr^T)$.

\begin{definition}\label{conv}
Let $f,g \in \src (\gr^T)$, we define a function $f*g$ in
$\gr^T$ by 
\[ (f*g)(x,\xi,0)=\int_{T_x\gr_x}f(x, \xi -\eta ,0)
g(x,\eta,0)
d\mu_x(\eta)  \]
and
\[ (f*g)(\gamma,t)=  \int_{\gr_{s(\gamma)}}f(\gamma \circ \delta^{-1},t)
g(\delta ,t) t^{-q}d\mu_{s(\gamma)}(\delta) \]
for $t\neq 0$.
\end{definition}

We can enounce our main result.

\begin{theorem}
$*$ defines an associative product on $\src (\gr^T)$.
\end{theorem}

\begin{proof}
Remember we are assuming for the moment the proposition
\ref{elproducto}. Let $f,g\in \src (\gr^T)$. We let $F:=(f,g)$ the
function in $(\gr^T)^{(2)}$ defined by:
$$(f,g)(x,\xi,\eta,0)=f(x,\xi,0)\cdot g(x,\eta,0)$$
and
$$(f,g)((\gamma,t),(\delta,t))=f(\gamma,t)\cdot g(\delta,t)$$
for $t\neq 0$. Now, from the Leibnitz formula for the derivate of a
product it is immediate that $(f,g)\in \src ((\gr^T)^{(2)})$. Finally,
by definition we have that
$$m_{r,c}((f,g))=f*g,$$
hence, thanks to proposition \ref{elproducto}, 
$f*g$ is a well defined element in $\src (\gr^T)$.

For the associativity of the product, let us remark that when one
restrict the product to $\ci_c (\gr^T)$, this coincides with the
product classically considered on $\ci_c (\gr^T)$ (which is
associative, see for example \cite{Pat}). The associativity
for $\src (\gr^T)$ is proved exactly in the same way that for $\ci_c (\gr^T)$.
\end{proof}

We have then to prove proposition \ref{elproducto}. We are going
to start locally. Let $U\in \Rr^p \times
\Rr^{q} \times \Rr^q$ be an open set and $V=U\cap\Rr^p \times \{0\}
\times \{0\}$. Let 
$P: \Rr^p \times \Rr^{q} \times \Rr^q \rightarrow \Rr^p \times \Rr^{q}$ the
canonical projection $(x,\eta, \xi) \mapsto (x,\eta)$. We set 
$U'=P(U)\in  \Rr^p \times \Rr^{q}$, then $U'$ is also an open subset,
$V\cong U'\cap \Rr^p \times \{0\}$ and $P|_{V}=Id_{V}$. 
We denote also by $P$ the restriction 
$P:U\rightarrow U'$. We have as in lemma \ref{ellema2} a $\ci$ map 
$\tilde{P}:\Omega_{V}^{U} \rightarrow \Omega_{V}^{U'}$,
which in this case is explicitly written by
\[ \tilde{P}(x,\eta,\xi,t)=(x,\eta,t)  \]
We define $\tilde{P}_{r,c}:\src (\Omega_{V}^{U}) \rightarrow 
\src (\Omega_{V}^{U'})$ as follows
\[ \tilde{P}_{r,c}(F)(x,\eta,t)
=\int_{ \{ \xi \in \Rr^q:(x,\eta,\xi,t)\in \Omega_{V}^{U} \} }
F(x,\eta,\xi,t)d\xi.  \]
Let us prove the following lemma.
\begin{lemma}\label{ptilde}
$\tilde{P}_{r,c}: \src (\Omega_{V}^{U})\rightarrow
\src (\Omega_{V}^{U'})$ is well defined.
\end{lemma}

\begin{proof}
The first observation is that the integral in the definition
of $\tilde{P}_{r,c}$ is always well defined. Indeed, we deduce it
from the next two points:
\begin{itemize}
\item[$\cdot$]For $t=0$, $\xi \mapsto F(x,\eta,\xi,0) \in \sw
  (\Rr^q)$.
\item[$\cdot$]For $t\neq 0$, $\xi \mapsto F(x,\eta,\xi,t) \in \ci_c (\Rr^q)$.
\end{itemize}
Once we can derivate under the integral symbol, we obtain that 
$\tilde{P}_{r,c}(F)\in \ci (\Omega_{V}^{U'})$. Then, 
we just have to show that $\tilde{P}_{r,c}(F)$ verifies the two
conditions of the definition
\ref{ladef}. For the first, if $K \subset U\times [0,1]$ is the 
compact conic support of $F$, then it is enough to put 
\[ K'=(P\times id_{[0,1]})(K)  \]
in order to obtain a conic compact subset of $U'\times [0,1]$ relative
to $V$ and to check that $K'$ is the compact conic support of 
$\tilde{P}_{r,c}(F)$ .
Let us now verify the condition $(s_1)$. Let $k,m\in \Nat$, $l\in \Nat^p$ 
and $\beta \in
\Nat^{q}$. We want to find $C_{(k,m,l,\beta)}>0$ such that
\[ (1+\| \eta \|^2)^k \| \partial_{x}^{l}\partial_{\eta}^{\beta}
\partial_{t}^{m}\tilde{P}_{r,c}(F)(x,\eta ,t)\| 
\leq C_{(k,m,l,\alpha)}   \]
For $k'\geq k+\frac{q}{2}$ and $\alpha=(0,\beta) \in \Rr^{q} \times \Rr^q$
we have by hypothesis that it exists $C'_{(k',m,l,\alpha)}>0$ such that
\[  \| \partial_{x}^{l}\partial_{\eta}^{\beta}
\partial_{t}^{m}F(x,\eta,\xi,t)\| 
\leq C' \frac{1}{(1+\| (\eta,\xi) \|^2)^{k'}}  \]
Then, we also have that
\[ \| \partial_{x}^{l}\partial_{\eta}^{\beta}
\partial_{t}^{m}\tilde{P}_{r,c}(F)(x,\eta ,t)\| 
\leq C' \int_{\{ \xi \in \Rr^q:(x,\eta,\xi,t)\in \Omega_{V}^{U} \}}
\frac{1}{(1+\| (\eta,\xi) \|^2)^{k'}}d\xi \]
\[ \leq C'\frac{1}{(1+\| \eta \|^2)^{\frac{q}{2}-k'}} 
\int_{ \{ \xi \in \Rr^q \} }
\frac{1}{(1+\| \xi \|^2)^{k'}}d\xi
 \leq C \frac{1}{(1+\| \eta \|^2)^{k}} \]
with
$$C=C'\cdot \int_{ \{ \xi \in \Rr^q \} }
\frac{1}{(1+\| \xi \|^2)^{k'}}d\xi .$$
\end{proof}

We can now give the proof of the proposition \ref{elproducto}.

\begin{proof}[Proof of \ref{elproducto}]
Let us first fix some notation. We suppose dim $\gr=p+q$ and
dim $\go=p$, in particular this implies that dim $\gr^{(2)}=p+q+q$. 
Let $(\Uo,\phi)$ and $(\Uo',\phi')$ be $\go$-slices in
$\gr^{(2)}$ and $\gr$ respectively such that the following diagram commutes
\[
\xymatrix{
\Uo \ar[r]^m \ar[d]_{\phi} & \Uo' \ar[d]^{\phi'} \\
U \ar[r]_P & U', }
\] 
where $P:\Rr^p \times \Rr^q \times \Rr^q \rightarrow \Rr^p \times \Rr^q$
is the canonical projection (as above) and $P(U)=U'$. This is possible
since $m$ is a surjective submersion.
Now, we apply the DNC construction to the diagram above to obtain,
thanks to the functoriality of the construction, the following 
commutative diagram
\[
\xymatrix{
\Dnc_{\Vo}^{\Uo} \ar[r]^{\Dnc(m)} \ar[d]_{\tilde{\phi}} & 
\Dnc_{\Vo'}^{\Uo'} \ar[d]^{\tilde{\phi'}} \\
\Dnc_{V}^{U} \ar[r]_{\Dnc(P)}  & 
\Dnc_{V'}^{U'}  \\
\Omega_{\phi} \ar[u]^{\Psi_{\phi}} \ar[r]_{\tilde{P}} &
\Omega_{\phi'} \ar[u]_{\Psi_{\phi'}},
}
\]
where $\Psi_{\phi}$ and $ \tilde{\phi}$ are as in section
\ref{dnc}. 
Let $g\in \src (\Dnc_{V}^{U})$,
we define
$$
P_{r,c}(g)(x,\eta ,t) = \left\{ 
\begin{array}{cc}
\int_{\Rr^q}g(x,\eta ,\xi,0)d\xi &\mbox{ if } t=0 \\
\\
\int_{\{ \xi \in \Rr^q : (x,\eta ,\xi,t) \in \Omega_{V}^{U}\}}
g(x,\eta ,\xi,t)t^{-q}d\xi &\mbox{ if } t\neq 0\\
\end{array}\right.
$$
Then, from the last commutative diagram, we get that 
$$P_{r,c}(g)= \tilde{P}_{r,c}(g\circ \Psi_{\phi})\circ
(\Psi_{\phi'})^{-1},$$ hence, thanks to \ref{ptilde}, 
we can conclude that we have a well defined linear map
\[ P_{r,c}: \src (\Dnc_{V}^{U}) \rightarrow \src (\Dnc_{V'}^{U'}).\]
We now use the proposition \ref{gtilde} to write 
\[ \src (\Dnc_{\Vo}^{\Uo})=\{ h\in \ci(\Dnc_{\Vo}^{\Uo}): 
h\circ \tilde{\phi}^{-1} \in \src (\Dnc_{V}^{U}) \}, \]
and so for $h\in \src (\Dnc_{\Vo}^{\Uo})$, we see that
\[ P_{r,c}(h\circ \tilde{\phi}^{-1})\circ \tilde{\phi'} 
\in \src (\Dnc_{\Vo'}^{\Uo'}).  \]
We use again the last commutative diagram to see that
\[ m_{r,c}(h)= P_{r,c}(h\circ \tilde{\phi}^{-1})\circ \tilde{\phi'}. \]
We then have a well defined linear map
\[ m_{r,c}: \src (\Dnc_{\Vo}^{\Uo}) \rightarrow \src
(\Dnc_{\Vo'}^{\Uo'}).  \]
To pass to the global case we only have to use the decomposition of 
$\src (\Dnc_{\go}^{\gr^{(2)}})$ and of 
$\src (\Dnc_{\go}^{\gr})$ as in (\ref{decomposicion}), and of course
 the invariance under diffeomorphisms (proposition
\ref{gtilde}).
\end{proof}

Let us recall that we have well defined evaluation morphisms as in (\ref{e0}) 
and (\ref{e1}). In the case of a the tangent groupoid they are by
definition morphisms of algebras. Hence, the algebra $\src(\gr^T)$ is 
field of algebras over the closed interval $[0,1]$, 
with associated fiber algebras,
\begin{center} 
$\sw(A\gr)$, at $t=0$, and
\end{center}
\begin{center}
$\ci_c(\gr)$ for $t\neq 0$.
\end{center}
It is very interesting to see what this means in the examples given in 
\ref{exgt}.

\section{Further developements}

Let $\gr \rightrightarrows \go$ be a Lie groupoid. In index theory for
Lie groupoids the tangent groupoid has been used to define the
analytic index associated to the group, as a morphism 
$K^0(A^*\gr) \rightarrow K_0(C_{r}^{*}(\gr))$ (see \cite{MP}) or as a
$KK$-element in \linebreak $KK(C_0(A^*\gr),C_{r}^{*}(\gr))$ (see \cite{HS}), 
this can be done because one has the following 
short exact sequence of $C^*$-algebras 
\begin{align}\label{se*}
0\rightarrow C_{r}^{*}(\gr \times (0,1]) \longrightarrow 
C_{r}^{*} (\gr^T)
\stackrel{e_0}{\longrightarrow} C_0(A^*\gr) \longrightarrow 0,
\end{align}
and because of the fact that the $K$-groups of the algebra $C_{r}^{*}(\gr
\times (0,1])$ vanish (homotopy invariance). The index defined at
the $C^*$-level has proven to be very useful (see for example \cite{CS}) 
but extracting numerical invariants from it, 
with the existent
tools, is very difficult. In non commutative geometry, and also in
classical geometry, the tools for obtain more explicit invariants are
more developed for the 'smooth objects'; in our case this means 
the convolution algebra $\ci_c(\gr)$, where we can for example apply
Chern-Weil-Connes theory. Hence, in some way, the indices defined in 
$K_0(\ci_c(\gr))$ are more refined objects and for some
cases it would be preferable to work with them. Unfortunately this
indices are not good enough, since for example they are not homotopy
invariants; in \cite{Co2} Alain Connes discusses this and also other
reasons why it is not enough to keep with the $\ci_c$-indices. 
The main reason to construct the algebra $\src ({\gr^T})$ is that it
gives an intermediate way between the $\ci_c$-level and the
$C^*$-level and will allow us in \cite{Ca} to define another analytic
index morphism associated to the groupoid, with the advantage that
this index will take values in a group that allows to do pairings with
cyclic cocycles and in general to apply Chern-Connes theory to it. The
way we are going to define our index is by obtaining first a short
exact sequence analogue to (\ref{se*}), that is, a sequence of the
following kind 
\begin{align}\label{se}
0\rightarrow J \longrightarrow \src (\gr^T)
\stackrel{e_0}{\longrightarrow} \sw (A^*\gr)) \longrightarrow 0.
\end{align}
The problem here will be that we do not dispose of the advantages of
the $K-$theory for $C^*$-algebras, since the algebras we are
considering are not of this type 
(we do not have for example homotopy invariance).

\section{References}

\renewcommand{\refname}{}    

\vspace*{-36pt}              

\end{document}